\documentclass[preprint]{ejpecp} 

\usepackage[T1]{fontenc}
\usepackage[utf8]{inputenc}

\usepackage{bbm}
\usepackage[retainorgcmds]{IEEEtrantools}

\usepackage{enumerate}  

\SHORTTITLE{A generalisation of the BDG inequalities}

\TITLE{A generalisation of the Burkholder-Davis-Gundy inequalities \
    \thanks{Saul Jacka gratefully acknowledges funding received from the Alan Turing Institute for their financial support under the 
EPSRC grant EP/N510129/1.}
\thanks{We are most grateful to an anonymous referee for very helpful suggestions and for a shortening of the proof of the main theorem.}} 

\AUTHORS{%
  Ma. Elena Hern\'andez-Hern\'andez\footnote{Department of Mathematics, Imperial College London, UK
    \BEMAIL{m.hernandez-hernandez@imperial.ac.uk}}
  \and 
  Saul D. Jacka\footnote{Alan Turing Institute and Department of Statistics, University of Warwick, UK \BEMAIL{s.d.jacka@warwick.ac.uk}
  }}


\AMSSUBJ{60G07; 60H05} 

\SUBMITTED{August 12, 2021} 
\ACCEPTED{October 9, 2022} 

\VOLUME{0}
\YEAR{2022}
\PAPERNUM{0}
\DOI{10.1214/YY-TN}

\ABSTRACT{Consider a \cad local martingale $M$ with square brackets $[M]$. In this paper, we provide upper and lower bounds for expectations of the type $\E \pq{M}^{q/2}_{\tau}$, for any stopping time $\tau$ and $q\ge 2$,  in terms of predictable processes.  This result can be thought of as  a  Burkholder-Davis-Gundy type inequality in the sense that it can be used to relate the expectation of the running maximum  $|M^*|^q$ to the expectation of the dual previsible projections of  the relevant powers of the associated jumps of $M$. The case for a class of moderate functions is also discussed.}

\numberwithin{theorem}{section}

\usepackage{color}

\usepackage{float}

\numberwithin{equation}{section}

\theoremstyle{remark}
\def\defto{\stackrel{\text{def}}{=}}

\newcommand{\E}{\mathbb{E}}

\def\half{\frac{1}{2}}
\def\qt{p}
\def\Dk{D^{(r)}}
\def\Df{D^{(2,F)}}

\def\Af{A^{(2,F)}}
\def\Atof{A^{(F)}}
\def\cf{c_F}
\def\Cf{C_F}
\def\tf{{\tilde F}}

\def\Zq{Z^{(q)}}
\def\Sq{S^{(q)}}

\newcommand{\cad}{\text{c\`adl\`ag }}

\newcommand{\lee}{\,\, \le \,\,}

\newcommand{\pb}[1]{\,\left (\, #1 \, \right )}
\newcommand{\pq}[1]{\,\left [\, #1 \, \right ] }
\newcommand{\ab}[1]{\left |\, #1 \, \right |}
\newcommand*\diff{\mathop{}\!\mathrm{d}}

\begin{document}


\section{Introduction}
In the context of stochastic calculus, the celebrated Burkholder-Davis-Gundy (BDG) inequalities play an important role in the estimation of moments of local martingales and, thus, more generally for  semimartingales and associated stochastic integrals. The BDG inequalities relate the maximum of a local martingale $M_t^* := \sup_{s\le t} |M_s|$ to its quadratic variation $[M]$.  More precisely, for any local martingale $M$ with $M_0 = 0$ and for any $1 \le q < \infty$, there exist universal  positive constants (independent of the local martingale $M$)  $c_q$ and $C_q$  such that, for any stopping time $\tau$,

\begin{equation}\label{BDG}
c_q\, \E\pb{ \pq{M}^{q/2}_{\tau}} \lee \E\pb{ \ab{M_{\tau}^*}^q} \lee C_q\, \E\pb{ \pq{M}^{q/2}_\tau}.
\end{equation} 
If $M$ is a continuous local martingale, the  inequalities in \eqref{BDG} hold for any $q> 0$ (see, e.g., \cite[p. 83]{LiptserS1989}, \cite[Theorem IV.74, p.226]{Protter}). 
{The generalisation for convex moderate functions $F$ is given by 
\begin{equation}\label{BDGF}
c_F \, \E\pb{  F \pb{ \pq{M}^{1/2}_{\tau}}}\,\, \lee \,\,\E\pb{ F \pb{M_{\tau}^*}} \,\, \lee \,\,C_F \, \E\pb{  F \pb{ \pq{M}^{1/2}_{\tau}}},
\end{equation}
for universal constants $c_F>0$ and $C_F>0$, see, e.g., \cite[Theorem 42.1, p.93]{RWilliams2}, \cite[Theorem 2.1]{Lenglart1980}.}

One drawback when applying  the BDG inequality \eqref{BDG} is  that, recalling that
\[
[M]\,\,=\,\,<M^c>\,\,+\sum\limits_{0<s\leq t}(\Delta M_s)^2 \quad \text{ and }\quad
\Delta[M]_t=(\Delta M_t)^2,
\]
 one may be dealing with a process where the compensator for the squared jumps is well understood but the jumps themselves are not. Moreover, the version of (\ref{BDG}) where $[M]$ is replaced by $<M>$, its dual previsible projection (or compensator) is, in general, false for $q>2$ when $M$ is a discontinuous local martingale, since, in general,
\[
\text{there exists } c_q:\; \E\bigl([M]_\tau^{\frac{q}{2}}\bigr)\leq c_q \E\bigl(<M>_\tau^{\frac{q}{2}}\bigr)\text{ only for }q\leq 2
\]
  (see \cite[Item (4.b'), Table 4.1, p.162]{BJY1986}).

In this note we prove the existence of universal positive constants $C_q$ and $c_q$ such that, for $q\geq 2$,

\[
c_q \, \E \pq{\max \left \{\,<M>_{
\tau}^{\frac{q}{2}},\,A_{\tau}^{(\frac{q}{2})} \,\right \}} \,\, \lee\,\,  \E\pb{[M]_\tau^{\frac{q}{2}}}\,\, \lee \,\,C_q \,  \E \pq{\max \left \{\,<M>_{\tau}^{\frac{q}{2}},\,A_{\tau}^{(\frac{q}{2})} \,\right \}},
\]
where  $<M>$ is the angle brackets of $M$, and $A^{(r)}$ is the dual previsible projection of the process 
\[ t \mapsto  \sum_{0<s\le t} |\Delta M_s|^{2r}, \quad \quad r \ge 1.
\]
Hence, our main results  compare the moments of the quadratic variation process of any local martingale with associated predictable processes (see Theorem \ref{BDGrev}). A generalisation for a class of moderate functions is also given in Theorem \ref{BDGmod}. \\
Since the BDG inequalities in \eqref{BDG} relate the running maximum of a local martingale to its quadratic variation, the further application of our results in Corollary \ref{BDGpred} allows us to derive estimates for the $q$th moments of the running maximum $M^*$ in terms of corresponding predictable processes. For the case $0<q < 2$, it is known that when $M$ is a locally square integrable local martingale then (see, e.g., \cite[Theorem 5, p. 69]{LiptserS1989}):
\[ 
  \E\pb{ \ab{M_{\tau}^*}^q} \, \lee \,\frac{4-q}{2-q} \,  \E \pq{<M>_{\tau}^{\frac{q}{2}}}, 
\]
 whereas, if $M$ is a continuous local martingale, then
 \begin{equation}\label{BDGq2C}
\frac{2-q}{4-q}\,\E \pq{<M>_{\tau}^{\frac{q}{2}}}\, \lee \, \E\pb{ \ab{M_{\tau}^*}^q} \, \lee \,\frac{4-q}{2-q} \,  \E \pq{<M>_{\tau}^{\frac{q}{2}}}. 
 \end{equation}
For any moderate function $F$, Lenglart \emph{et.al.}  \cite[Section 2, p. 37]{Lenglart1980} provide estimates for discontinuous local martingales  with jumps bounded by a locally bounded predictable increasing process $D$, i.e. $|\Delta M| \le D$:
\[ 
\E\pb{ F \pb{M_{\tau}^*}} \,\, \lee \,\,C \, \E\pb{  F \pb{ <M>^{\half}_{\tau} + D_{\tau}}},
\] 
 and 
 \begin{equation}\nonumber
 \E\pb{  F \pb{ <M>^{\half}_{\tau}}} \lee \,\, c\, \E\pb{ F \pb{M_{\tau}^* + D_{\tau}}}.
\end{equation}
A special class of local martingales that are encompassed by our results are those  obtained as stochastic integrals with respect to continuous local martingales such as Brownian motion, or with respect to local martingales given by  compensated Poisson random measures. In these two cases the representation of the dual previsible projections $A^{\left (\frac{q}{2}\right)}$ can be written explicitly and one can thus derive (for the one-dimensional case) the well-known estimates for the running maximum of associated local martingales (see, e.g., Kunita (2004) \cite{Kunita2004-InBook}, Applebaum (2009) \cite[Theorem 4.4.22 and Theorem 4.4.23]{a}, Marinelli and R\"ockner (2014)   \cite{Marinelli2014}). 
Extensions and generalisations to the multidimensional case and Hilbert-space-valued case, and applications of our results to the estimates of $q$th moments of semimartingales will be provided in a forthcoming paper. 

\section{Main Result}
 Consider a filtered complete probability space $(\Omega, \mathcal{F}, \mathbb{F}, \mathbb{P})$, on which all our processes are defined.    $\mathcal{A}^+_{loc}$ denotes the class  of locally integrable, non-decreasing, \cad and adapted processes.  Let us also recall that if $A \in \mathcal{A}^+_{loc}$, then there is a predictable process $A^p \in \mathcal{A}^+_{loc}$, called the \emph{dual previsible projection} (or compensator) of $A$, which is unique up to an evanescent set, and which is characterised by making $A - A^p$ a local martingale (or equivalently   $\E (A_{T\wedge T_n}^p) = \E (A_{T\wedge T_n})$ for all stopping times $T$ and for a localising sequence $T_1,T_2,\ldots$),  see, e.g., \cite[Theorem I.3.17, p. 32]{JacodS1987}.

The main result is the following

\begin{theorem}\label{BDGrev}
For any \cad local martingale $M$ and for  $r\ge 1$,
define the adapted, increasing process $D^{(r)}$ by
\[
D^{(r)}_{t} \defto \sum_{0<s\le t} |\Delta M_s|^{2r},
\]
and define  $A^{(r)}$ to be the dual previsible projection of $\Dk$ (whenever it exists). For any $q \geq 2$
define the process
\[
t\,\,\mapsto S_t^{(q)}(M) \defto \max
 \left \{<M>_t^{\frac{q}{2}},\,\,A_t^{(\frac{q}{2})}\right\}.
\]

There exist universal constants $c_q>0$ and $C_q >0$ such that for all stopping times $\tau$ and local martingales $M$,

\begin{equation}\label{E:T}
 c_q \,\,\E \pb{\Sq_\tau (M)} \,\,\, \lee \,\, \, \E \pb{[M]_\tau^{\frac{q}{2}}}\,\,\,\lee \,\,\,C_q\,\,\E\pb{\Sq_\tau (M)}.
\end{equation}
\end{theorem}
\begin{remark}
Note that
\[
<M>\,\,=\,\,<M^c>\,+\,A^{(1)},
\]
so that, in particular,
\[
S^{(2)}=\,\,<M>.
\]
\end{remark}
\begin{proof}[Proof of Theorem \ref{BDGrev}]
Let us first observe that, since $M$ is a \cad local martingale, it has a countable number of jumps and, thus, using the fact that the $\ell_q$ spaces satisfy $\ell_2 \subset \ell_{2r}$ for all $r\ge 1$, we obtain: 
\begin{equation}\nonumber
\pb{ \sum_{0<s\le t}  |\Delta M_s|^{2r}}^{\frac{1}{2r}} \le   \pb{ \sum_{0<s\le t}  |\Delta M_s|^{2}}^{\half},\quad \text{for all }\,\,\, r\ge 1,
\end{equation}
which in turn implies that
 \begin{equation}\nonumber
 D^{(r)}_t \le \left (\, D_t^{(1)}\,\right)^r, \quad \quad \text{ for all }\,\,\, r \ge 1.
 \end{equation}
In particular, for $r = q/2$ with $q\ge 2$, we have that
\begin{equation}\label{E:ellM}
\E \pb{  \sum_{0<s\le t}  |\Delta M_s|^q }  \,\,\, \le \,\,\, \E \pb{ \pb{  \sum_{0<s\le t}  \pb{\Delta M_s}^2  }^{\frac{q}{2}}}  \,\, \leq \,\,  \E \pb{  [M]_{t}^{\frac{q}{2}} }.
\end{equation}
To prove the left-hand inequality in \eqref{E:T}, assume that  $\E \pb{  [M]_{\tau}^{\frac{q}{2}} }  \,<\, \infty$ (as otherwise there is nothing to prove).  Thus, the inequality  \eqref{E:ellM} implies that  $t\mapsto D^{(\frac{q}{2})}_{\tau \wedge t} \in \mathcal{A}^+_{loc}$  and so its dual predictable projection $A^{(\frac{q}{2})}_{\tau \wedge \cdot}$ exists (\cite[Theorem I.3.17, p.32]{JacodS1987}). Now recall the standard result that if $D$ is an increasing, adapted and locally integrable process started at 0 and $A$ is its dual previsible projection,  then for any $p\geq 1$
\[
\E[A_\tau^p]\lee p^p\,\E[D^p_\tau],
\]
(see, e.g., \cite[Theorem 4.1]{Lenglart1980}). 
Applying this inequality to the processes  
$D^{(p)}$ and $A^{(p)}$, and using the fact that 
\[
[M]^p\,\, \geq \,\, \max \,\left \{ D^{(p)},\,<M^c>^p \right\},
\]
 the left-hand inequality in \eqref{E:T} follows.

For the right-hand inequality, first observe that, since $S^{(2)}\, = \, <M>$, the result is trivial for $q=2$; so suppose that $q>2$. We assume that $\E\pb{\Sq_\tau(M)}< \infty$, as otherwise there is nothing to prove. 
Now note that, using the fact that (for $x,y\geq 0$) $x\vee y\leq x+y\leq 2(x\vee y)$ and observing that 
$[M]=[M^d]+<M^c>$, it is enough to prove the right-hand inequality for purely discontinuous martingales.

To simplify notation we denote $\frac{q}{2}$ by $\qt$.
Define the process $\Zq$ by
\begin{IEEEeqnarray*}{ll}
\Zq_{}=[M]^{\qt}.
\end{IEEEeqnarray*}
 
Notice that $f:x\mapsto|x|^\qt$ is $C^{1}$ so, using the fact that
$|x+y|^p\lee k_p(|x|^p+|y|^p)$ where $k_p := \max\{2^{p-1},1\}$, for any  $p>0$,   the Mean Value Theorem implies that there exists $\theta \in (0,1)$ such that

\begin{IEEEeqnarray}{ll}\nonumber
\Delta \Zq_{t}&=f([M]_{t})-f([M]_{t-})\\\nonumber
&=p\biggl([M]_{t-}+\theta (\Delta M_{t})^2\biggr)^{p-1}(\Delta M_{t})^{2}\\
& \le \,\,p\, k_{p-1}\,\,\pb{[M]_{t-}^{p-1}(\Delta M_{t})^{2} + (\Delta M_{t})^{2p}}. \label{E:Zq} 
\end{IEEEeqnarray}

Now $\Zq$ is increasing from 0 and increases only by jumps, so the estimate \eqref{E:Zq} implies that
\begin{IEEEeqnarray}{ll}\nonumber
\Zq_t
&\lee \int_{0+}^t b_p\,[M]^{p-1}_{s-}\diff [M]_s\,+\,b_p\,D^{(p)}_{t}, \label{Z:02q}
\end{IEEEeqnarray}
with $b_p:= p\, k_{p-1}$. 
 
We take dual previsible projections  and evaluate at $\tau$ to obtain
\begin{IEEEeqnarray}{ll}\label{ineq2}
\E\pb{[M]_\tau^p}&\lee b_p\,\E \pb{\int_0^\tau [M]^{p-1}_{s-}\diff <M>_s+\,\,A^{(p)}_{\tau}}\nonumber\\
&\lee b_p\E\pb{ [M]^{p-1}_\tau<M>_\tau+\,\,A^{(p)}_{\tau}}\nonumber\\\,\,
&\lee b_p\, \pb{||[M]_\tau||_p^{p-1}\,\,||<M>_\tau||_p\,+\E (\,A^{(p)}_\tau\,)}.
\end{IEEEeqnarray}
Dividing both sides of \eqref{ineq2} by $\E[\Sq_\tau]$ we see that, setting $z= \pb{\frac{\E\pb{[M]_\tau^p}}{\E(\Sq_\tau)}}^{\frac{1}{p}}$,
\[
z^p\lee b_p(z^{p-1}+1)\defto g(z).
\]
Denoting by $c_q$ the largest root of $g(x)=x^p$, we see that 
the right-hand inequality in  \eqref{E:T}  holds, as required.

\end{proof}

The following example is well-known to show the failure of $<M>$ to control the moments of $[M]$ when $M$ is a discontinuous local martingale. 
\begin{example}\label{ex1}
Take $M$ to be a compound Poisson process of unit rate with jump times $T_1,T_2,\ldots$, and with jump-sizes $X_1,X_2,\ldots$, where
\[
\E[X_1]=0,\quad \; \E[X_1^2]=1\quad \text{ and } \quad \E[X_1^4]=\infty.
\]
It follows that $M$ is a square-integrable martingale with 
\[
[M]_t=\sum_n1_{\{T_n\leq t\}}X^2_n,
\]
and with
\[
<M>_t=t.
\]
Clearly
\begin{equation}\label{big}
\E\pb{[M]_t^2}=\infty,
\end{equation}
so there is no $c$ such that
\[
\E\pb{[M]_T^2}\,\,\leq \,\,c\,\E \pb{<M>_T^2}.
\]
However, 
$\E[D^{(2)}_t]=\infty$, for any $t>0$, as implied by \eqref{E:T} and \eqref{big}.
\end{example}

The standard BDG inequality and Theorem \ref{BDGrev} imply the following
\begin{corollary}\label{BDGpred}
For any \cad local martingale $M$ and $q\ge 2$, there are universal constants $c_q>0$ and $C_q >0$ such that for all stopping times $\tau$ 
\begin{equation}\label{BDGp}
c_q \, \E \pq{\max \left \{\,<M>_{
\tau}^{\frac{q}{2}},\,A_{\tau}^{(\frac{q}{2})} \,\right \}} \, \lee\,  \E\pb{ \ab{M_{\tau}^*}^q} \, \lee \,C_q \,  \E \pq{\max \left \{\,<M>_{\tau}^{\frac{q}{2}},\,A_{\tau}^{(\frac{q}{2})} \,\right \}}. \end{equation} 
\end{corollary}
\begin{remark}
Note that when $M$ is a continuous local martingale then $A^{(\frac{q}{2})} \equiv 0$ and one recovers the inequalities in  \eqref{BDGq2C} for any $q \ge 2$.
\end{remark}

\section{The case of moderate functions}
A function $F : \mathbb{R}_+ \to \mathbb{R}_+$ is said to be moderate if it is continuous and increasing, $F(x) = 0$ and if, for some  (and then for every)  $\alpha>1$, the following  growth condition holds\footnote{Equivalently, if the condition $\sup_{x>0} \frac{F(\alpha x)}{F(x)}  \,\,<\,\,\, \infty $ holds.}: 
\[\text{for some }c>0,\;\;F (\alpha x) \lee c\,F(x)\text{ for all }x>0.
\]
If $F$ is convex with right derivative $f$, then a necessary and sufficient condition for $F$ to be moderate is that  
\begin{equation}\label{E:expF}
q := \sup_{x>0} \frac{x f(x)}{F(x)} \,\,<\,\, \infty,
\end{equation}
see, e.g., \cite[Section 1]{Lenglart1980}. 
Here, $q$ is known as the exponent of $F$. 

If \eqref{E:expF} holds, then for all $\alpha > 1$, 
\begin{equation}\label{E:GthF}
 \sup_{x>0} \frac{F(\alpha x)}{F(x)} \lee \alpha^q.\end{equation}
 In particular, 
 the power function $x \mapsto x^p$ for $p\ge1$ is a moderate convex function and, in this case, $c$ can be taken as $\alpha^p$ and the exponent of $F$ is equal to $p$.

One is naturally led to ask whether there is a generalisation of Theorem \ref{BDGrev}, to some class of  moderate  functions, of the form
\begin{IEEEeqnarray}{ll}\nonumber\label{moderate3}
\cf\E\left [\max \left \{F\pb{<M>^\half_\tau},\Af_\tau \right \} \right ]&\lee \E \left [F \pb{[M]^\half_\tau}\,\right ]\\\nonumber
&\lee \Cf \E \left [\max \left \{F\pb{<M>^\half_\tau},\Af_\tau \right \}\right ],\\
\end{IEEEeqnarray}
where $A^{(2,F)}$ is the dual previsible projection of the corresponding process  $D^{(2,F)}$ (defined in \eqref{moderate1} below).

If one takes $F: x \mapsto x^q$, for any $q \ge 2$, (these are, of course,  convex moderate functions), then Theorem \ref{BDGrev}  guarantees that the two-sided inequality \eqref{moderate3} holds for any continuous local martingale $M$. However, we can see that for discontinuous local martingales \eqref{moderate3} may not hold for $q <2$. 
 Indeed, a small change to Example \ref{ex1} shows that the left hand side of \eqref{moderate3} cannot hold in general: 
\begin{example}\label{ex2}
Take $M^n$ to be a compound Poisson process of unit rate with jump times $T_1,T_2,\ldots$, and with jump-sizes $X_1,X_2,\ldots$, where
\[
\E[X_1]=0,\; \,\,\quad \E[|X_1|]=1\,\,\quad \text{ and }\,\,\quad \E[X_1^2]=n.
\]
It follows that $M^n$ is a square-integrable martingale with 
\[
[M^n]_t=\sum_n1_{\{T_n\leq t\}}X^2_n, \quad \quad \; \;\E \pb{[M^n]_t^\half}\,\,\leq \,\,t
\]
and with
\[
<M^n>_t\,\,=\,\,nt.
\]
Clearly, taking $F$ to be the identity, it follows that, for the left-hand inequality in \eqref{moderate3} to hold, we would need $c_F\geq n$ for any $n$, contradicting the finiteness of $c_F$.
\end{example}

\begin{theorem}\label{BDGmod}
Suppose that $F$ is a strictly increasing and convex moderate function. Define $\Af$ and $\Atof$ to be the dual previsible projections of 
\begin{equation}\label{moderate1}
\Df \defto \sum\limits_s F(|\Delta M_s|) \quad \text{ and } \quad D^{(F)}\defto \sum\limits_s F \pb{(\Delta M_s)^2}, 
\end{equation} respectively. There are universal constants $\cf>0$ and $\Cf>0$ such that

\begin{equation}\label{moderate2} \E \left [\,F \pb{[M]^\half_\tau}\,\right ]\lee \,\,\Cf\, \E \left [ \,\max \left \{F \pb{<M>^\half_\tau},\,\Af_\tau \right \}\, \right ].
\end{equation}
and
\begin{equation}\label{moderate4}
\cf\,\E\left [\,\,\max \left \{F(<M>_\tau),\,\Atof_\tau \right \}\, \right ]\,\,\lee\,\, \E \left [\,F([M]_\tau)\,\right ]\,\,
\end{equation}
\end{theorem}
\begin{proof} Let us first assume that $F$ is $C^1$. We use, without further comment, the facts that $F'(\cdot)$ is increasing, $[M]_\cdot$ and $<M>_\cdot$ are increasing, $\Delta[M]_t^\half\leq |\Delta M_t|$ and $<M^c>\;\leq [M]$.

Proceeding as in the proof of Theorem \ref{BDGrev}, we obtain first
\[
V_t\defto \sum_{0<s\le t} \Delta F\left ([M]^\half_s \right )\lee  \sum_{0<s\le t}  F'\left ([M]^\half_{s} \right )\Delta [M]^\half_s\lee  \sum_{0<s\le t}  F'\pb{[M]^\half_{s}}|\Delta M_s|.
\]
By \eqref{E:expF}, for all $x$, there exists $q>0$ such that $F'(x) \lee q \frac{F(x)}{x}$, and so
\begin{equation}\nonumber
V_t\lee \sum_{0<s\le t}  F'\left ([M]^\half_{s}\right)\pb{|\Delta M_s|}\lee q\,    \sum_{0<s\le t}  \frac{F\left ([M]^\half_{s} \right )}{[M]^\half_s}|\Delta M_s|.
\end{equation}
The convexity of $F$ together with \eqref{E:expF} then implies that
\begin{IEEEeqnarray}{ll} \nonumber 
V_t&\lee  \frac{q}{2} \sum_{0<s\le t}  \frac{F\left (2[M]^\half_{s-}\right)+F(2|\Delta M_s|)}{[M]^\half_s}|\Delta M_s|\\ \nonumber
&\lee  b_F\, \sum_{0<s\le t}  \frac{F\left ([M]^\half_{s-}\right )\,+\,F\pb{|\Delta M_s|}}{[M]^\half_s}|\Delta M_s|\\ \nonumber
&\lee b_F\, \sum_{0<s\le t}  \left ({F'([M]^\half_{s-} )}|\Delta M_s|+F\pb{|\Delta M_s|}\right ),
\end{IEEEeqnarray}
 where $b_F := q\, 2^{q-1}$. 
 
Now,  \cite[Theorem II.31, p.78]{Protter} implies 
\[
F\left([M]^\half_t\right)=\int_{0+}^tF'\left([M]^\half_{s-}\right)\frac{\diff <M^c>_s}{[M]^\half_s}\,\,+\,\,V_t,
\]
so, evaluating at $\tau$, using the previous inequality and then taking dual previsible projections yield
\begin{IEEEeqnarray}{ll}\label{F3}
\E\left [F\left ([M]^\half_\tau \right )\right ]&\lee b_F\,\E \pq{\int_{0+}^\tau F'\left ([M]^\half_{t-} \right )\frac{\diff <M^c>_t}{[M]_t^\half}\,+\, \Af_\tau}\nonumber\\
&\lee b_F\,\E \pq{F'([M]_{\tau})\int_0^\tau\frac{\diff <M^c>_t}{<M^c>_t^\half}+\Af_\tau}\nonumber\\
&\lee 2b_F\,\E \pq{F'\left ([M]^\half_{\tau} \right)<M^c>^\half_t+\,\Af_\tau}.
\end{IEEEeqnarray}
Denote the convex conjugate of $F$ by $\tf$, then the generalisation of Young's inequality implies %

\[
xy\lee \tf \pb{\frac{x}{\mu}}\,+\,F(\mu y)
\]
for any $ x,y,\mu>0$.
We also have (see \cite[Lemma 1.1.1]{BJY1986}) 
\[
\tf \pb{\frac{F'(x)}{2}}\lee F(x).
\]
Thus, taking  $\mu=\max \{ 2\,b_F,1\}$, we see that
\begin{IEEEeqnarray}{ll} \nonumber
\E\left [F' \left ([M]^\half_{\tau} \right)<M^c>^\half_\tau \right ]&\lee \E \pq{\tilde{F} \pb{\frac{F'\left ([M]^\half_{\tau}\right )}{2\mu}}\,+\,F\left (2\mu <M^c>^\half_\tau \right)}\\ \nonumber
&\lee \frac{1}{\mu}\E \left [F\left ([M]^\half_\tau \right)+(2\mu)^qF\left (<M^c>^\half_\tau \right) \right ].
\end{IEEEeqnarray}
Substituting this into the inequality \eqref{F3} we obtain
inequality \eqref{moderate2}.

To obtain the  inequality \eqref{moderate4}, observe that we need only show that, for some $c>0$,
\[
\E[F(<M>_\tau)]\lee c\,\E[F([M]_\tau)],
\]
since $\Delta F([M]_t)\geq F(\Delta[M]_t)$ by convexity of $F$.
Now, for some $\theta_{\cdot} \in (0,1)$,
\begin{IEEEeqnarray}{ll}
F(<M>_t)&\,=\, \int_{0+}^t F'(<M>_s)\diff <M>^c_s\,+\, \sum_{0<s\le t}  F'(<M>_{s-}\,+\,\theta_s \Delta<M>_s)\Delta<M>_s\nonumber\\
&\lee \int_{0+}^t F'(<M>_s) \diff <M>^c_s+ \sum_{0<s\le t}  F'(<M>_{s})\Delta<M>_s \nonumber \\
&=\int_{0+}^t F'(<M>_s) \diff <M>_s.
\end{IEEEeqnarray}
Stopping at $\tau$ and then taking expectations  we obtain that, for any $\mu>0$,
\begin{IEEEeqnarray*}{ll}
\E[F(<M>_\tau)]&\lee \E \pq{\int_{0+}^\tau F'(<M>_t)\diff [M]_t}\lee  \E[F'(<M>_\tau)[M]_\tau]\\
&\lee   \E \pq{\tf \pb{\frac{F'(<M>_\tau)}{2\mu}}+F(2\mu[M]_\tau)}\\
&\lee  \E \pq{\frac{1}{\mu}F(<M>_\tau)+(2\mu)^qF([M]_\tau)}. 
\end{IEEEeqnarray*}
Taking $\mu=2$ we obtain the desired result.

For the general case, we can now proceed by approximation. Take a sequence of $C^2$ convex functions $\{F_n\}_{n\ge 0}$ such that $F_n  \uparrow F$, and with the corresponding derivatives increasing to the left derivative $F'_{-}$ of $F$. The sequence $\{F_n\}_{n\ge 0}$ can be constructed in the standard way as a convolution \lq\lq on the left'' of $F$ with an appropriate sequence of scaled versions of a (positive) mollifier $\phi$ (i.e. a  $C^{\infty}$ function with compact support that integrates to 1):
\[
F_n(x)=\int_0^{\infty}F(x-t){\phi(nt)}ndt,
\] 
where we define $F(y)=0$ for $y<0$.
Then the previous case implies that for each $F_n$, the inequalities in \eqref{moderate2} {and \eqref{moderate4}} hold and so the result follows by letting $n\to \infty$ and using the  Monotone Convergence Theorem.
\end{proof}

\begin{remark}
By taking $F$ appropriately in Theorem \ref{BDGmod}, we recover the inequalities in Theorem \ref{BDGrev}. Indeed, \eqref{E:T} follows from \eqref{moderate4} and \eqref{moderate2} by taking $F: x \mapsto x^{q/2}$ and $F: x \mapsto x^{q}$, for $q>2$, respectively.  It is worth observing that Theorem \ref{BDGmod} also extends the RHS inequality of \eqref{E:T} in Theorem \ref{BDGrev} to the case $q \in (1,2)$.   
\end{remark}
The BDG inequality in \eqref{BDGF} and Theorem \ref{BDGmod} imply the following
\begin{corollary}\label{BDGpredF} Let $F$ be a strictly increasing and convex moderate function. 
For any \cad local martingale $M$, there exist universal constants $c_F,\,\,C_F >0$ such that for all stopping times $\tau$ 
\begin{IEEEeqnarray}{rl}\label{BDGp}
  \E\pb{ F \pb{M_{\tau}^*}} \, &\leq \Cf\, \E \left [ \max \left \{F \pb{<M>^\half_\tau},\Af_\tau \right \} \right ]\\
  \E \left [ \,\max \left \{F \pb{<M>_\tau},\,\Atof_\tau \right \}\, \right ]&\leq \,\,\cf\,  \E\pb{ F \pb{(M_{\tau}^*)^2}},
   \end{IEEEeqnarray} 
   so that (\ref{moderate3}) does hold {\em if} $F:x\mapsto G(x^2)$, where $G$ is a strictly increasing, convex, moderate function.
\end{corollary}


\bibliographystyle{amsplain}




\end{document}